\documentclass{article}\begin{document}\centerline{\bf Some Notes on the Solutions of non Homogeneous Differential Equations} 
\centerline{\bf with Constat Coefficients}\vskip .4in

\centerline{\bf Nikolaos Bagis} 
\centerline{\bf Department of Informatics}
\centerline{\bf Aristotele University of Thessaloniki}
\centerline{\bf Greece}
\centerline{\bf nikosbagis@hotmail.gr}

\abstract{We solve some forms of non homogeneous differential equations in one and two dimensions. By expanding the solution into whell-posed closed form-Eisenstein series the solution itself is quite simple and elementary. Also we consider Fourier series solutions of linear differential operator equations. In the third section we study operators which are functions of the Leibnitz derivative. The last result is the complete solution of a non homogenus 2-degree ODE with linear coeficients. The non homogenous part is an arbirtary function of $L_2(\bf R\rm)$}.

\[
\]

\section{The Divisor Sums and ODE}

\textbf{Proposition 1.} If $x$, is positive real number and $f$ is analytic in (-1,1), with $f(0)=0$, then
\begin{equation}
\exp\left(\int^{x}_{\infty}f(e^{-t})dt\right)=\prod^{\infty}_{n=1}(1-e^{-nx})^{\frac{1}{n}\sum_{d|n}\frac{f^{(d)}(0)}{d!}\mu(n/d)}
\end{equation}
Where $\mu$ is the Moebius function. See also and [Ap] chapter 2 for the Moebius as also for other multiplicative functions.\\ 
\textbf{Proof.}
See [B]
\[
\]
Examples of such identities are\\
1) Let $\frac{f^{(n)}(0)}{n!}=1[n]=1$ if $n=1$, 0 else. $X(n)=\frac{1}{n}\mu(n)$
$$\prod^{\infty}_{n=1}(1-q^n)^{\frac{\mu(n)}{n}}=e^{-q}$$
2) Let $\frac{f^{(n)}(0)}{n!}=n$, $n=1,2,\ldots$ then $f(x)=\frac{x}{(x-1)^2}$ and\\ $X(n)=\frac{1}{n}\sum_{d|n}d\mu(n/d)=\frac{\phi(n)}{n}$
$$\prod^{\infty}_{n=1}(1-q^n)^{\frac{\phi(n)}{n}}=e^{\frac{q}{q-1}}$$
Where $\phi(n)$ is the Euler Totient function.\\   
\textbf{Proposition 2.} If $A(n)$ is arbitrary sequence of numbers we have for $x>0$
\begin{equation}
\frac{d^{\nu}}{dx^{\nu}}\left(\sum^{\infty}_{n=1}\frac{\sum_{d|n}A(d)\mu(n/d)}{e^{nx}-1}\right)=\sum^{\infty}_{n=1}\frac{\sum_{d|n}A(d)(-d)^{\nu}\mu(n/d)}{e^{nx}-1}
\end{equation}
\textbf{Proof.}\\
See also [B]
\[
\]
We will use Proposition 1 to find the solution of the $N-th$ degree linear differential equation
\begin{equation}
\sum^{N}_{\nu=0}a_{\nu}\frac{d^{\nu}}{dx^{\nu}}u(x)=\sum^{\infty}_{n=1}\frac{\sum_{d|n}C(d)\mu(n/d)}{e^{nx}-1}=\sum^{\infty}_{n=1}C(n)e^{-nx}
\end{equation}
\[
\]
\textbf{Lemma 1.} Set 
\begin{equation}
P(x)=\sum^{N}_{\nu=0}a_{\nu}x^{\nu}
\end{equation}
then the solution of (3) is
\begin{equation}
u(x)=\sum^{\infty}_{n=1}\frac{\sum_{d|n}\frac{C(d)}{P(-d)}\mu(n/d)}{e^{nx}-1}=\sum^{\infty}_{n=1}\frac{C(n)}{P(-n)}e^{-n x} 
\end{equation}
\textbf{Proof 1.}
From Proposition 2, it is clear that if $$u(x)=\sum^{\infty}_{n=1}\frac{\sum_{d|n}A(d)\mu(n/d)}{e^{nx}-1}
$$   
for a certain $A(k)$, then (3) becomes
$$\sum^{N}_{\nu=0}a_{\nu}(x)\left(\sum^{\infty}_{n=1}\frac{\sum^{\infty}_{d|n}A(d)\mu(n/d)}{e^{nx}-1}\right)^{(\nu)}=\sum^{\infty}_{n=1}\frac{\sum_{d|n}C(d)\mu(n/d)}{e^{nx}-1}$$
or
$$\sum^{N}_{\nu=0}a_{\nu}(x)\sum^{\infty}_{n=1}\frac{\sum_{d|n}A(d)(-d)^{\nu}\mu(n/d)}{e^{nx}-1}=\sum^{\infty}_{n=1}\frac{\sum_{d|n}C(d)\mu(n/d)}{e^{nx}-1}$$
or
$$\sum^{\infty}_{n=1}\frac{\sum_{d|n}A(d)P(-d)\mu(n/d)}{e^{nx}-1}=\sum^{\infty}_{n=1}\frac{\sum_{d|n}C(d)\mu(n/d)}{e^{nx}-1}$$
Hence it must be
\begin{equation}
\sum_{d|n}A(d)P(-d)\mu(n/d)=\sum_{d|n}C(d)\mu(n/d)
\end{equation}
Eq.6 shows clearly that $$A(n)=C(n)/P(-n)$$
Also one can see that we have 
$$
\sum^{\infty}_{n=1}\frac{W(n)}{e^{nx}-1}=\sum^{\infty}_{n=1}\left(\sum_{d|n}W(d)\right)e^{-nx}
$$
and
$$
\sum^{\infty}_{n=1}\frac{\sum_{d|n}W(d)\mu(n/d)}{e^{nx}-1}=\sum^{\infty}_{n=1}W(n)e^{-nx}
$$
\textbf{Proof 2.}
Let $$u(x)=\sum^{\infty}_{n=1}u_ne^{-nx}$$ then setting into (3) the above expansion we get the same result in a more easy way. 
\[
\]
\textbf{Theorem 1.} If
$$
\sum^{N}_{\nu=0}a_{\nu}u^{(\nu)}u(x)=\sum^{\infty}_{n=1}\frac{C(n)}{e^{nx}-1}
$$
Then
$$u(x)=\sum^{\infty}_{n=1}\frac{\sum_{d|n}C(d)}{P_x(n)}e^{-nx}$$, where $$P_x(w)=\sum^{N}_{k=0}a_k w^k$$
\textbf{Proof.}
As in Lemma 1.
\[
\]
Next we proceed with the 2-dimension problem with a similar way. 
\[
\]
We set 
\begin{equation}
P_x(w)=\sum^{N}_{k=0}a_k w^k
\end{equation} 
and 
\begin{equation}
P_y(w)=\sum^{M}_{l=0}b_l w^l
\end{equation}
Also
\begin{equation}
G(x,y)=\sum^{\infty}_{k,m=1}\frac{c(k,m)}{(e^{kx}-1)(e^{my}-1)}
\end{equation}
\textbf{Theorem 2.} The equation 
\begin{equation}
\sum^{N,M}_{k,l=0}a_{k}b_{l}u^{(k),(l)}(x,y)=G(x,y)
\end{equation}
have solution
\begin{equation}
u(x,y)=\sum^{\infty}_{n,m=1}\frac{S(n,m)}{(e^{nx}-1)(e^{my}-1)}\end{equation}
where
\begin{equation}
S(n,m)=\sum_{d|n,\delta|m}B(d,\delta)\mu(n/d)\mu(m/\delta)
\end{equation}
and
\begin{equation}
B(n,m)=\frac{1}{P_x(-n)P_y(-m)}\sum_{k|n,r|m}c(k,r)
\end{equation}
\textbf{Proof.}
Let $$
u(x,y)=\sum^{\infty}_{n=1}\frac{\sum_{d|n}A_d(y)\mu(n/d)}{e^{nx}-1}
$$
Then differentiating with respect to $x$ we get 
$$
u^{(k),(0)}(x,y)=\sum^{\infty}_{n=1}\frac{\sum_{d|n}A_d(y)(-d)^{k}\mu(n/d)}{e^{nx}-1} \eqno{(a)}
$$
then with respect to $y$ we get
$$
u^{(k),(l)}(x,y)=\sum^{\infty}_{n=1}\frac{\sum_{d|n}A_d^{(l)}(y)(-d)^{k}\mu(n/d)}{e^{nx}-1} \eqno{(b)}
$$
but
$$
A_d^{(l)}(y)=\sum^{\infty}_{m=1}\frac{\sum_{\delta|m}B(d,\delta)(-\delta)^{l}\mu(m/\delta)}{e^{my}-1} \eqno{(c)}
$$
combining the above we get the result.
\[
\]
\textbf{Note.}
The polynomials $P$ apparently must have no solutions in natural numbers. 
\[
\]
\textbf{Examples.}\\
\textbf{1)} Set 
$$P_x(w)=1+\sqrt{2}w+w^2$$
also set
$$P_y(w)=1+\sqrt{2}w$$
Then the equation 
$$
u(x,y)+\sqrt{2}u_{x}(x,y)+u_{xx}(x,y)+\sqrt{2} u_{y}(x,y)+2u_{xy}(x,y)+\sqrt{2} u_{xxy}(x,y)=$$
$$=\sum^{\infty}_{n=1}\frac{1}{e^{nx}-1}\sum^{\infty}_{m=1}\frac{1}{e^{my}-1}
$$ 
have solution   
$$u(x,y)=\left(\sum^{\infty}_{n=1}\frac{\sigma_0(n)e^{-nx}}{1-\sqrt{2}n+n^2}\right)\left(\sum^{\infty}_{m=1}\frac{\sigma_0(m)e^{-my}}{1-\sqrt{2}m}\right)$$
where $$\sigma_{\nu}(n)=\sum_{d|n}d^{\nu}$$
\textbf{2)} If
$$P_x(w)=2+\sqrt{3}w+w^3$$
$$P_y(w)=1+\sqrt{3}w+w^2$$
then the solution of
$$2u+\sqrt{3}u_{x}+u_{xxx}+2\sqrt{3}u_{y}+3u_{xy}+\sqrt{3}u_{xxxy}+2u_{yy}+\sqrt{3}u_{xyy}+u_{xxxyy}=$$
$$=\sum^{\infty}_{n,m=1}\frac{\log(n+m)}{(e^{nx}-1)(e^{my}-1)}$$
is
$$
u(x,y)=\sum^{\infty}_{n,m=1}\frac{\sum_{d|n}\sum_{\delta|m}B(d,\delta)\mu(n/d)\mu(m/\delta)}{(e^{nx}-1)(e^{my}-1)}
$$
where $$B(n,m)=\frac{\sum_{d|n}\sum_{\delta|m}\log(d+\delta)}{(2- \sqrt{3}n-n^3)(1-\sqrt{3}m+m^2)}$$
Observe that in this example we are not able to split the solution into two parts in $x$ and $y$.\\
\textbf{Note.} For no confusion the form of the equation is defined by $$\sum^{N,M}_{n,m=0}a_nb_mu^{(n),(m)}(x,y)$$ where the $a_n$ and $b_n$ are respectively that of $P_x$ and $P_y$. 

\section{Series Solutions}

If $F$ is an operator such that
$$F(x)=\sum^{\infty}_{k=0}\frac{(t\lambda)^k}{k!}\frac{d^k}{dx^k}=e^{\lambda t\frac{d}{dx}}\eqno{(d)}$$
We will try to solve the equation $$\frac{\partial}{\partial t}u(x,t)=\frac{\partial}{\partial x}u(x,t)$$

Assume that 
$$
u(x,t)=e^{t\frac{\partial}{\partial x}}f(x)
$$
then $$\frac{\partial}{\partial t}u(x,t)=e^{t\frac{\partial}{\partial x}}\frac{\partial}{\partial x}f(x)=\frac{\partial}{\partial x}e^{t\frac{\partial}{\partial x}}f(x)=\frac{\partial}{\partial x}u(x,t)
$$
Hence the operator $\exp\left(t\frac{\partial}{\partial x}\right)$ produces the solution.
From Eq.(3) and Lemma one can take the limit $N\rightarrow \infty$ then $$
F\left(\frac{d}{dx}\right)u(x,t)=\exp\left(-t\frac{\partial}{\partial x}\right)u(x,t)=f(x)=\sum^{\infty}_{n=1}C(n)e^{-nx}
$$ 
thus according to Lemma the solution of the above equation must be $$u(x,t)=\sum^{\infty}_{n=1}\frac{C(n)}{e^{t n}}e^{-n x}$$     
Using the parameter $\lambda$ which can take and complex values one can arrive to the conclusion that 
\begin{equation} 
\frac{\partial^{\nu}}{\partial t^{\nu}}u(x,t)+\frac{\partial^{\nu}}{\partial x^{\nu}}u(x,t)=0
\end{equation}
have solution 
\begin{equation}
u(x,t)=\sum^{\infty}_{n=1}\frac{C(n)}{\exp\left[e^{i\pi/\nu} n t\right]}e^{-nx}
\end{equation}
One can see that, a solution of Schrodingers equation
\begin{equation}
u_t(x,t)=-\frac{\partial^m}{\partial x^m}u(x,t)+V(t)u(x,t)
\end{equation}
is
\begin{equation}
u(x,t)=\sum^{\infty}_{n=1}\frac{C(n)e^{-nx}}{\exp\left[(-n)^m t-f(t)\right]}
\end{equation}
where 
\begin{equation}
f(x)=\int^{x}_{c}V(w)dw
\end{equation}
This is the general case in which the potential depends only in time.
\[
\] 
Let again
\begin{equation}
P(x):=\sum^{N}_{k=0}a_kx^k
\end{equation}
Consider now the equation 
\begin{equation}
\sum^{N}_{k=0}a_ky^{(k)}(x)=P\left(\frac{d}{dx}\right)y(x)=x-\pi
\end{equation}
the solution is
\begin{equation}
y(x)=\frac{xa_0-a_0\pi-a_1}{a^2_0}+\sum^{N}_{k=1}C(k)e^{x\rho_k}
\end{equation}
Where $\rho_k$ is the roots of $P(x)=0$ 
The same equation have solution according to the Theorems of section 1:
\begin{equation} 
y(x)=\sum_{n\in \bf Z^{*}\rm}\frac{i}{n P(in)}e^{inx}
\end{equation}
\[
\]
An interesting question is how one can extract from (21) and (22) the roots $\rho_k$.\\
Anyway when if we let $N\rightarrow \infty$, then\\
\textbf{Theorem 3.}\\
\textbf{i)} 
\begin{equation} 
y(x)=\sum_{n\in \bf Z^{*}\rm}\frac{i}{n F(in)}e^{inx}
\end{equation}
\textbf{ii)}
\begin{equation}
F\left(\frac{d}{dx}\right)y(x)=x-\pi
\end{equation} 
\textbf{iii)}
\begin{equation}
y(x)=\frac{xa_0-a_0\pi-a_1}{a_0^2}+\sum^{\infty}_{k=1}C(k)e^{\rho_k x}
\end{equation}
the $\rho_k$ are roots of Eq. $F(x)=x-\pi$. The function $F$ must have not integer roots in the imaginery line.
\[
\]
Now consider the function $F(x)=e^{-2\pi i x}+x-\pi-1$. It is $a_0=F(0)=-\pi$ and $a_1=-2\pi i+1$ and $\rho_k=k$. Hence the two representations are $$y(x)=\sum_{n\in \bf Z^{*}\rm}\frac{ie^{inx}}{n(e^{2n\pi}+in-\pi-1)}$$
$$y(x)=\frac{-x\pi+\pi^2+2\pi i-1}{\pi^2}+\sum^{\infty}_{n=-\infty}C(n)e^{inx}$$
the differential equation is $$y(x-2\pi i)+y'(x)-y(x)(\pi+1)=x-\pi
$$   
\textbf{Examples}\\  
\textbf{1)} $$F(x)=e^x+x+1$$ then 
$$y(x)=\sum_{n\in \bf Z^{*}\rm}\frac{i}{n(-n+1+e^{n})}e^{inx}$$
and also
$$y(x)=\frac{-1-\pi+x}{2}+e^{\pi ix}\sum^{\infty}_{k=-\infty}C(k)e^{2k\pi ix}$$ 
in this example the diferential eq. is
$$y(x+1)+y'(x)+y(x)=x-\pi$$
\textbf{2)} For $F(x)=\cosh(x)+x+1$ 
the DE is $$\frac{y(x+1)+y(x-1)}{2}+y'(x)+y(x)=x-\pi$$
with solution $$y(x)=\sum_{n\in \bf Z^{*}\rm}\frac{i}{nF(in)}e^{inx}$$
The series in a first view can not become more fast convergent. If we consider for example $$F(x)=\cos(x)+x+1$$ then $$y_{M}(x)=\sum_{|n|\leq M}\frac{i}{n(\cosh(n)+in+1)}e^{inx}$$
But for the diferential equation holds
$$\frac{y_{M}(x+i)+y_{M}(x-i)}{2}+y_{M}'(x)+y_{M}(x)-x+\pi=O\left(\frac{1}{M}\right)$$
If we try with $F(x)=\cos(\pi x)+x+1$ then we have 
$$y_{M}(x)=\sum_{|n|\leq M}\frac{i}{n(\cosh(n\pi)+in\pi+1)}e^{inx}$$
and for the diferential equation holds
$$\frac{y_{M}(x+i\pi)+y_{M}(x-i\pi)}{2}+y_{M}'(x)+y_{M}(x)-x+\pi=O\left(\frac{1}{M}\right)$$
which is the same.\\
Thus we can say that even the equations are not easy to solve numericaly, the solutions itself under certain conditions may behave very good i.e $y_{M}(x)$ is very fast convergent.\\ 
\textbf{3)} Let us consider now a curius case. The $L_2(\bf R\rm)$ function $h(x)=e^{x-e^{x}}$ then if $F(x)=h(x)+x+1$ (we dont need the roots), we have 
\begin{equation}
\sum^{\infty}_{l=0}\frac{(-1)^l}{l!}{y(x+l+1)}+y'(x)+y(x)=x-\pi
\end{equation} 
and 
\begin{equation}
y(x)=\sum_{n\in \bf Z^{*}\rm}\frac{i}{n(in+1+e^{in-e^{in}})}e^{inx}
\end{equation}   
is the solution indeed. 

\section{Functions of $\frac{d}{dx}$}

We proceed with the following\\ 
\textbf{Lemma 2.} Let $$f(x)=\sum^{\infty}_{n=0}f_nx^n$$ be analytic function in $\bf R\rm$ such that for every $a,b>0$ there exist constant depending from $f$, $M_f$:
\begin{equation}
|f(x)|\leq M_f (1+|x|^a) e^{-b|x|}
\end{equation}
Let also $\phi(x)$ real valued function with values in $\bf C\rm$, such that for every $c>0$ there exist costant $M_{\phi}$:
\begin{equation}
|\phi(x)|\leq M_{\phi}|x|^c
\end{equation}  
then 
\begin{equation}
\int^{\infty}_{0}f(x)\phi(x)e^{-xs}dx=\sum^{\infty}_{k=0}f_{2k}\frac{\partial^{2k} (L\phi)(s)}{\partial s^{2k}}-\sum^{\infty}_{k=0}f_{2k+1}\frac{\partial^{2k+1} (L\phi)(s)}{\partial s^{2k+1}}
\end{equation}
\textbf{Proof.}\\ See [Ba]
\[
\]
If happens $L\phi=y(x)$ then 
\begin{equation} 
\int^{\infty}_{0}f(x)(L^{(-1)}y)(x)e^{-xs}dx=
$$
$$=\sum^{\infty}_{k=0}\frac{f^{(2k)}(0)}{(2k)!}\frac{d^{2k}y(s)}{ds^{2k}}-\sum^{\infty}_{k=0}\frac{f^{(2k+1)}(0)}{(2k+1)!}\frac{d^{2k+1}y(s)}{ds^{2k+1}} 
\end{equation}
hence we can write
$$ 
\int^{\infty}_{0}f(x)(L^{(-1)}y)(x)e^{-xs}dx=f_e\left(\frac{d}{ds}\right)y(s)-f_o\left(\frac{d}{ds}\right)y(s)$$
\begin{equation}
\int^{\infty}_{0}f(-x)(L^{(-1)}y)(x)e^{-xs}dx=f\left(\frac{d}{ds}\right)y(s)
\end{equation} 
\textbf{Theorem 4.}\\i) It holds   
\begin{equation}
f\left(\frac{d}{ds}\right)y(s)=L\left(f(-x)(L^{(-1)}y)(x)\right)(s)
\end{equation}
ii) The solution of 
\begin{equation}
f\left(\frac{d}{ds}\right)y(s)=g(s)
\end{equation}
is
\begin{equation}
y(s)=\int^{\infty}_{0}\frac{(L^{(-1)}g)(x)}{f(-x)}e^{-x s}dx
\end{equation}   
\[
\]
This theorem shows clearly that we can find solutions in integral-closed-form, of the general not homogeneous equation (if existing the Laplace transforms).  
\[
\]
\textbf{Examples.}\\
\textbf{1)} If hapens $g(x)=1/x^2$, then $(L^{(-1)}g)(x)=x$ and thus
the solution of 
\begin{equation}
f\left(\frac{d}{ds}\right)y(s)=1/s^2
\end{equation}
is
\begin{equation}
y(x)=\int^{\infty}_{0}\frac{x}{f(-x)}e^{-xs}dx
\end{equation}
where the form of the equations (36) and (37) is that of (34) and (35).  
In the same way as in the above example we can set other values for $g(x)$\\
\textbf{2)}
Set $$h(x)=\sum^{\infty}_{k=0}c_kx^{\nu_k}e^{-l_kx}$$
Then the equation
\begin{equation}
\sum^{\infty}_{k=0}c_ky^{(\nu_k)}(x-l_k)=g(x)
\end{equation}
have solution 
\begin{equation}
y(x)=\int^{\infty}_{0}\frac{(L^{(-1)}g)(w)}{h(-w)}e^{-w x}dw
\end{equation}   
This method is like solving (38) with Fourier or Laplace transforms but we avoid some restrictions of $y$ and $g$ to be in $L_2(\bf R\rm)$. Note also that it is solved with Laplace theory.\\
\textbf{3)} We try now to evaluate $T=\frac{1}{1+\frac{d}{dx}}$. 
Let $Ty(x)=g(x)$, then $$\frac{1}{1+\frac{d}{dx}}=1-\frac{d}{dx}+\frac{d^2}{dx^2}-\frac{d^3}{dx^3}+\ldots$$     
$$g(x)=Ty(x)=\frac{1}{1+\frac{d}{dx}}y(x)=1-\frac{dy}{dx}+\frac{d^2y}{dx^2}-\frac{d^3y}{dx^3}+\ldots$$  
we use (39) and get $$y(x)=\int^{\infty}_{0}(L^{(-1)}g)(w)(1-w)e^{-x w}dw$$
solving with respect to $g$ we have
$$ 
g(w)=Ty(w)=\int^{\infty}_{0}\frac{(L^{(-1)}y)(x)}{1-x}e^{-xw}dx
$$
In general holds
\begin{equation}
h\left(\frac{d}{dx}\right)f(x)=\int^{\infty}_{0}(L^{(-1)}f)(w) h(-w)e^{-xw}dw
\end{equation}  
under certain conditions of convergence. For example if $y(x)=$polynomial in $x$.  
Also for
$h(x)=\log(1+x)$, then
$$\log\left(1+\frac{d}{dx}\right)y(x)=-\sum^{\infty}_{n=1}\frac{(-1)^n}{n}\frac{d^n}{dx^n}y(x)=$$
$$=\int^{\infty}_{0}(L^{(-1)}y)(w)\log(1-w)e^{-xw}dw$$ 
\[
\]
Thus for example if we consider the equation $$f\left(\frac{d}{dx}\right)y(x)=g(x) \eqno{(e)}$$
with $f(x)=x+\log(1+x)$, then the differential equation $(e)$ is actualy the $$\frac{dy(x)}{dx}+\int^{\infty}_{0}(L^{(-1)}y)(w)\log(1-w)e^{-x w}dw=g(x)$$ and have solution $$y(x)=\int^{\infty}_{0}\frac{(L^{(-1)}g)(w)}{-w+\log(1-w)}e^{-wx}dw$$
\[
\]
Relation (40) is very useful if one can set a one to one relation between $h$ and a function of $y$. For example if one take $h(x)=e^x$ then for all the functions $$Q(x)+e^x$$ with $Q(x)=\sum^{N}_{k=0}a_kx^k$, the equation will be $$\sum^{N}_{k=0}a_k\frac{d^k}{dx^k}y(x)+y(x+1)=g(x)$$ will have the same solution type $$y(x)=\int^{\infty}_{0}\frac{(L^{(-1)}g)(w)}{Q(-w)+e^{-w}}e^{-xw}dw$$  
Can we say that for a given function $h(x)$, ($e^x$ in the examples) exists a unique form of $y(x)$ (such as $y(x+1)$ in the examples)? The exponential and polynomial functions behave very good, but what happens with other values of $h$. For example: 
Exist $h_2(x)$ giving us the form $y(x)^2$ in the differential equation?       
Then all the differential equations, with this term, will be solvable in integral forms with knowing one function only. But this seems not to happen. It happens with the differential functions, i.e.   $h\left(\frac{d}{dx}\right)y(x)=\log(1-\frac{d}{dx})$.\\
For example if we consider the equation 
$$2\frac{d^2}{dx^2}y(x)-a\frac{d}{dx}y(x)-\log\left(1-\frac{d}{dx}\right)y(x)=g(x)$$ the solution is  $$y(x)=\int^{\infty}_{0}\frac{(L^{(-1)}g)(w)}{2w^2+aw-\log(1+w)}e^{-xw}dw .$$
Another related equation is
$$y(x+i)+\log\left(1-\frac{d}{dx}\right)y(x)=g(x)$$
which have a solution $$y(x)=\int^{\infty}_{0}\frac{(L^{(-1)}g)(w)}{e^{-iw}+\log(1+w)}e^{-wx}dw$$
\[
\]
From the above examples and (39) and (40) one can see that the inversion with respect to some $g(x)$ is\\
\textbf{Theorem 5. (Inversion)}
\begin{equation}
h\left(\frac{d}{dx}\right)^{(-1)}=\frac{1}{h\left(\frac{d}{dx}\right)}
\end{equation}
This means if   
\begin{equation}
y(x)=\int^{\infty}_{0}\frac{(L^{(-1)}g)(w)}{h(-w)}e^{-xw}dw
\end{equation}
then
\begin{equation}
y^{(-1)}(x)=g(x)=\int^{\infty}_{0}(L^{(-1)}y)(w)h(-w)e^{-xw}dw
\end{equation}
Where $$h\left(\frac{d}{dx}\right)y(x)=g(x)=y_1(x)$$ and $$\frac{1}{h\left(\frac{d}{dx}\right)}y_1(x)=y(x)$$
\[
\] 
\textbf{Example.}\\ If $h(x)=e^{-x}+1$ then $h^{(-1)}(x)=\frac{1}{h(x)}=\frac{1}{e^{-x}+1}$. This means that if $g(x)=1/x^2$
$$y(x+1)+y(x)=g(x)$$ then $$y(x)=\frac{1}{4}\left(-\psi\left(1,1+\frac{x}{2}\right)+\psi\left(1,\frac{x+1}{2}\right)\right)$$
and $$g(x)=y_1(x)=\int^{\infty}_{0}(L^{(-1)}y)(w)(e^w+1)e^{-xw}dw=1/x^2$$
Where $\psi$ is the Polygamma function i.e $$\psi(z)=\frac{1}{\Gamma(z)}\frac{d\Gamma(z)}{dz}$$ and $$\psi(n,z)=\psi^{(n)}(z)$$ (see and Mathematica notes). 
The above example is trivial and can be solved with Laplace theory.
\[
\]
Now we will find a way to solve the equation
\begin{equation}
(a_1x+b_1)f''(x)+(a_2x+b_2)f'(x)+(a_3x+b_3)f(x)=g(x)
\end{equation}
where $f$, $g\in L_2(\bf R\rm)$.\\
Let the Fourier Transform of a function of $L_2(\bf R\rm)$ is $$\widehat{f}(\gamma)=\int^{\infty}_{-\infty}f(t)e^{-i t\gamma}dx$$
the Inverse Fourier Transform is $$f(x)=\frac{1}{2\pi}\int^{\infty}_{-\infty}\widehat{f}(\gamma)e^{i \gamma x}d\gamma$$
\textbf{Lemma 3.}
\begin{equation}
\int^{\infty}_{-\infty}f(x)x^ne^{-ix\gamma}dx=i^n(\widehat{f})^{(n)}(\gamma) .
\end{equation}
\begin{equation}
\widehat{(f^{(n)})}(\gamma)=(i\gamma)^n\widehat{f}(\gamma) .
\end{equation}
$$
\int^{\infty}_{-\infty}f'(x)A(x)e^{-i x\gamma}dx=
$$
\begin{equation}
=\int^{\infty}_{-\infty}f(x)A'(x)e^{-ix \gamma}dx+(-i\gamma)\int^{\infty}_{-\infty}f(x)A(x)e^{-ix\gamma}dx .
\end{equation}
$$
\int^{\infty}_{-\infty}f''(x)A(x)e^{-ix\gamma}=\int^{\infty}_{-\infty}f(x)A''(x)e^{-ix\gamma}dx+
$$
\begin{equation}
2(-i\gamma)\int^{\infty}_{-\infty}f(x)A'(x)e^{-ix\gamma}dx+(-i\gamma)^2\int^{\infty}_{-\infty}f(x)A(x)e^{-ix\gamma}dx .
\end{equation}
\textbf{Proof.}\\
The proof of (45) and (46) can be found in [Pa]. The relations (47) and (48) are obtained with integration by parts.\\   
\textbf{Theorem 6.}\\
When $f$, $g\in L_2(\bf R\rm)$ and $\lim_{|x|\rightarrow \infty}|f(x)x^{2+\epsilon}|=0$, $\epsilon>0$, equation (44) can reduced in into
\begin{equation}
(-ia_1\gamma^2+a_2\gamma+ia_3)\frac{\widehat{f}(\gamma)}{d\gamma}+(-b_1\gamma^2-2ia_1\gamma+ib_2\gamma+a_2+b_3)\widehat{f}(\gamma)=\widehat{g}(\gamma)
\end{equation} 
which is solvable.\\
\textbf{Proof.}\\
Take the Fourier Transform in both sides of (44) and use Lemma 3.

\[
\]

\newpage

\centerline{\bf References}

[Ap] T. Apostol. ''Introduction to Analytic Number Theory''. Springer-Verlag, New York, Berlin, Heidelberg, Tokyo 1976, 1984
 
[B] Nikos Bagis. ''Some results on the Theory of Infinite Series and Divisor Sums''. arXiv:0912.4815
 
[Ba] Nikos Bagis. ''Brief Research Notes on Transformation of Series and Special Functions''. arXiv:0907.1091 

[Pa]: A. Papoulis. ''The Fourier Integral and its Applications''. McGraw-Hill Publications. New York., 1962  

[S]: Murray R. Spiegel. ''Schaum's Outline of Theory and Problems of Fourier Analysis with Applications to Boundary Value Problems''. McGraw-Hill, Inc. 1974.

\end{document}